\documentclass{amsart}
\usepackage{amssymb}
\usepackage{amscd}
\usepackage{amsfonts}
\usepackage{latexsym}
\usepackage[all]{xy}
\usepackage{verbatim}

\newtheorem{theorem}{Theorem}[section]
\newtheorem{lemma}[theorem]{Lemma}
\newtheorem{proposition}[theorem]{Proposition}
\newtheorem{corollary}[theorem]{Corollary}
\theoremstyle{definition}
\newtheorem{definition}[theorem]{Definition}
\newtheorem{example}[theorem]{Example}

\newtheorem{remark}[theorem]{Remark}
\newtheorem{note}[theorem]{Note}
\newcommand{\id}{\text{id}}

\newcommand{\FPdim}{\text{FPdim}}

\newcommand{\Fun}{\text{Fun}}

\renewcommand{\Vec}{\text{Vec}}

\newcommand{\Hom}{\text{Hom}}

\newcommand{\Aut}{\text{Aut}}

\newcommand{\Rep}{\text{Rep}}

\newcommand{\tr}{\text{tr}\,}

\newcommand{\C}{\mathcal{C}}
\newcommand{\D}{\mathcal{D}}
\newcommand{\E}{\mathcal{E}}

\newcommand{\Z}{\mathcal{Z}}

\newcommand{\ot}{\otimes}

\begin{document}
\title[]
{Nilpotent fusion categories}

\author{Shlomo Gelaki}
\address{Department of Mathematics, Technion-Israel Institute of
Technology, Haifa 32000, Israel}
\email{gelaki@math.technion.ac.il}

\author{Dmitri Nikshych}
\address{Department of Mathematics and Statistics,
University of New Hampshire,  Durham, NH 03824, USA}
\email{nikshych@math.unh.edu}

\date{October 23, 2006}
\begin{abstract}
In this paper we extend categorically the notion of a finite 
nilpotent group to fusion categories. To this end, we first 
analyze the trivial component of the universal
grading of a fusion category $\C$, and then introduce the upper central series of $\C$.
For fusion categories with commutative Grothendieck rings (e.g., braided fusion categories) we also
introduce the lower central series.   
We study arithmetic and structural properties of nilpotent 
fusion categories, and apply our theory to modular categories and to 
semisimple Hopf algebras. In particular, we show that in the modular 
case the two central series are centralizers of each other in the sense of M.~M\"uger.
\end{abstract}
\maketitle

\begin{center}
{\bf Dedicated to Leonid  Vainerman on the occasion of his 60-th birthday }
\end{center}

\section{introduction}

The theory of fusion categories arises in many areas of
mathematics such as representation theory, quantum groups,
operator algebras and topology. 
The representation categories of 
semisimple (quasi-) Hopf algebras are important examples of fusion categories.
Fusion categories have been studied
extensively in the literature, and there exist many 
results concerning their structure and classification (see \cite{ENO} and
references therein). However, there are still many fundamental
open questions about fusion categories, which are motivated by
group theory and the theory of semisimple Hopf algebras, and it is desirable
to continue to develop the theory of fusion categories along these lines.

The purpose of this paper is to extend categorically the notion of a finite 
nilpotent group to fusion categories. For this end, we analyze the trivial component of the universal
grading of a fusion category $\C$, and then introduce the upper central series of $\C$. 
For fusion categories with commutative Grothendieck rings (e.g., braided fusion categories) we 
also introduce the lower central series.
We study arithmetic and structural properties of nilpotent
fusion categories, and apply our theory to modular categories and to
semisimple Hopf algebras. In particular, we show that in the modular
case the two central series are centralizers of each other in the sense of M.~M\"uger \cite{Mu}.


The organization of the paper is as follows.

In Section 2 we recall necessary definitions and results about based rings and
modules, fusion categories, and Frobenius-Perron dimensions,
which are needed in the sequel.

In Section 3, we first define the notion of the {\em adjoint}
subring $R_{ad}$ of a based ring $R$.  Then we prove in
Theorem \ref{constructing a group} that $R$ is naturally graded by a group 
$U(R)$
(which we call the {\em universal grading group} of $R$) and that $R_{ad}$ 
is the
trivial component of this grading. 
The adjective ``universal" above is justified
in Corollary \ref{faithful} where it is proved that any faithful grading of
$R$ arises from a quotient of $U(R)$. In case $R$ is the
Grothendieck ring of a fusion category $\C$ we prove that the
character group of the maximal abelian quotient of $U(R)$ is
isomorphic to the group of tensor automorphisms of the identity
functor of $\C$ (see Proposition \ref{maxab}). When $\C=\Rep(H)$,
the representation category of 
a semisimple Hopf algebra $H$, we have $\C_{ad}=\Rep(H/HK^+)$ where $K=\Fun(U(\C))$ is the maximal
central Hopf subalgebra of $H$ (see Theorem \ref{hopfad}).

Next we consider based rings of integer Frobenius-Perron
dimension. We show in Theorem~\ref{int} that any such based ring is
graded by an elementary abelian $2$-group, and that 
the objects of  integer dimension form the trivial component
of this grading.
This in particular implies that if  the  Frobenius-Perron dimension 
of a fusion category $\C$ is an odd integer then
the dimension of any object is an integer,
and hence $\C$ is equivalent to the  representation category  of
a  semisimple quasi-Hopf algebra (see Corollary
\ref{divprop}).

In Section 4 we use the construction of the adjoint based subring and adjoint
fusion subcategory to define the notions of the {\em upper central
series}, {\em nilpotency} and {\em nilpotency class} for based
rings and fusion categories, generalizing the corresponding notions 
for groups (see Definitions \ref{upper central series}, \ref{upp central series},
\ref{nilpotent}). 
We say that a semisimple Hopf algebra $H$ is nilpotent if the fusion category
$\Rep(H)$ of representations of $H$ is nilpotent. 
When $H$ is the group Hopf algebra of a finite group  this agrees with the classical
notion of a nilpotent group. Other examples of nilpotent fusion categories are
pointed fusion categories, Tambara-Yamagami categories, and $p$-fusion categories
(fusion categories of dimension $p^n$, $p$ a prime). It follows from \cite{ENO} 
that a nilpotent fusion category is pseudounitary and admits a canonical spherical structure.

Next, we define the notion of the {\em commutator} of a based subring of a
{\em commutative} based ring and of a fusion subcategory of a fusion category 
with commutative Grothendieck ring (see Definition \ref{com} and Definition \ref{comfus}). We
then use it to define 
the {\em lower central series} of {\em commutative} 
based rings and fusion categories with commutative Grothendieck ring (e.g., braided
fusion categories), generalizing the corresponding notion for groups (see Definitions 
\ref{lower central series}, \ref{low central series}).
Finally, in Theorem \ref{upp=low} we compare the two series and prove that the upper
one converges to $\mathbb{Z}1$ at step $n$ if and only if the lower one converges to $R$ 
at step $n$. This generalizes a classical result in group theory, cf.\  \cite{H}.

In Section 5 we consider indecomposable $R$-modules and prove that they 
decompose as $R_{ad}$-modules with index set being a transitive $U(R)$-set
(see Proposition \ref{constructing a Aset}). We use this to prove that in a
nilpotent  fusion category $\C$ 
the square of the Frobenius-Perron dimension of any simple object divides 
the Frobenius-Perron dimension of $\C_{ad}$, see Theorem \ref{inductive step} and Corollary 
\ref{nilpotent divisibility} (this generalizes a well-known property of nilpotent groups).
When applied  to semisimple Hopf algebras of a prime power dimension,  Corollary 
\ref{nilpotent divisibility} extends a result in \cite{MW}.

In Section 6 we focus on modular categories. We prove in 
Theorem \ref{Schneider's iso} that for a 
modular category $\C$ its universal grading group $U(\C)$ is canonically
isomorphic to the character group of the group of invertible objects
of  $\C$. For a pseudounitary modular category $\C$ we  find out 
in Theorem \ref{low = centralizer up} that the lower and upper series 
of $\C$ are related via the operation of taking centralizers, as defined by M.~M\"uger
in \cite{Mu}. As consequences we obtain that $\C_{ad}$ is the centralizer of 
the maximal pointed subcategory of $\C$, and that $\C$ contains a ``large" 
symmetric fusion subcategory (Corollaries \ref{cor1}, \ref{cor2}). Finally, we 
prove in Theorem \ref{brnil} that a braided fusion category $\C$ is nilpotent if and only
if its center $\mathcal{Z}(\C)$ is nilpotent. Moreover, if the nilpotency
class of $\C$ is $c$ then the nilpotency class of $\mathcal{Z}(\C)$ is at
most $2c$.

\par
\noindent {\bf Acknowledgments.} We are grateful to Pavel Etingof for sharing with us
his insights about the universal grading of a fusion category and for many helpful comments. 
We also thank Deepak Naidu, Victor Ostrik, and Leonid Vainerman for useful discussions.
The first author is grateful to the Departments of Mathematics
at the University of New Hampshire and MIT for their warm hospitality 
during his Sabbatical. His research was supported by the BSF grant no. 2002040.
The research of the second author was supported by the NSF grant DMS-0200202.

\section{Preliminaries}

Throughout the paper we work with an algebraically closed field $k$ of characteristic
$0$. All categories considered in this paper are finite, abelian, semisimple, and 
$k$-linear. All rings and modules have finite $\mathbb{Z}$-rank.

\subsection{Based rings and modules}
Let $\mathbb{Z}_+$ be the semi-ring of non-negative integers. We
will recall some definitions from \cite{O}. Let $R$ be a ring
with identity which is a finite rank $\mathbb{Z}$-module. A
{\em $\mathbb{Z}_+$-basis} of $R$ is a basis $B =\{X_i\}_{i\in
\mathcal{I}}$ such that $X_iX_j = \sum_{k\in \mathcal{I}} c_{ij}^k\, X_k,$ where
$c_{ij}^k\in \mathbb{Z}_+$. An element of $B$ will be called {\em
basic}.

Let us define a non-degenerate symmetric $\mathbb{Z}$-valued inner
product on $R$ as follows. For all elements $X= \sum_{i\in
\mathcal{I}}\,a_i X_i$ and $Y= \sum_{i\in \mathcal{I}}\,b_i X_i$
of $R$ we set
\begin{equation}
\label{pairing}
(X, Y) =\sum_{i\in \mathcal{I}}\, a_ib_i.
\end{equation}

\begin{definition}
\label{based ring} A {\em based ring} (or {\em multi-fusion ring})
is a pair $(R,\,B)$ consisting of a ring $R$ with a
$\mathbb{Z}_+$-basis $B =\{X_i\}_{i\in \mathcal{I}}$ satisfying
the following properties:
\begin{enumerate}
\item[(1)] There exists a subset $\mathcal{I}_0\subset \mathcal{I}$ such that $1 =\sum_{i\in \mathcal{I}_0}\,
X_i$.
\item[(2)] There is an involution $i \mapsto i^*$ of $\mathcal{I}$ such that
the induced map $X= \sum_{i\in \mathcal{I}} a_iX_i \mapsto X^*=
\sum_{i\in \mathcal{I}}a_i X_{i^*}$ satisfies $$(XY,\,Z) = (X,
ZY^*) = (Y,\, X^*Z)$$ for all $X,Y,Z\in R$.
\end{enumerate}
\end{definition}

In what follows we will usually denote a based ring $(R, B)$
simply by $R$ assuming that some  basis $B=\{X_i\}_{i\in
\mathcal{I}}$ is fixed.

\begin{remark}
It follows from Definition~\ref{based ring} that $X\mapsto X^*$ is
an anti-automorphism of the ring $R$. Also, for each $i\in
\mathcal{I}$ there is a unique $i_0\in \mathcal{I}_0$ such that
the coefficient of $X_{i_0}$ in the linear decomposition of the
product $X_i X_{i^*}$ is equal to $1$ while  
the coefficient of $X_j$ is $0$  for
each $j\in \mathcal{I}_0$ such that $ j\neq i_0$.
\end{remark}

A based ring is called {\em unital} (or {\em fusion}) if
$|\mathcal{I}_0|=1$.

\begin{note}
In this paper all based rings will be assumed unital. All our definitions
and results can be easily extended to the non-unital case.  
\end{note}

By a {\em based subring} of a based ring $(R,\,B)$ we will always understand a
based ring
$(S, C)$ where $C$ is a subset  of $B$ and $S$ is a subring of $R$.

A  {\em based (left) module} over a based ring $R$ is a free
$\mathbb{Z}$-module $M$ endowed with a fixed basis $\{V_j\}_{j\in
\mathcal{J}}$ such that $X_iV_j =\sum_{k\in J}\, d_{ij}^k V_k,\, d_{ij}^k
\in \mathbb{Z}_+$, and $d_{ij}^k = d_{i^*k}^j$ for all $i\in
\mathcal{I}$ and $j\in \mathcal{J}$.  The last condition simply means that
for all $X\in R$ the $\mathbb{Z}_+$-matrices of multiplication by
$X$ and $X^*$ in the basis $\{V_j\}_{j\in \mathcal{J}}$ are transposes of
each other. We will call an element of the basis of a based module {\em basic}.

Similarly to \eqref{pairing}, for all elements 
$V =\sum_{j\in \mathcal{J}}\, m_j V_j$ and
$U=  \sum_{j\in \mathcal{J}}\, n_j V_j$ of $M$ we set
\begin{equation}
(V, U) =\sum_{i\in \mathcal{J}}\, m_jn_j.
\end{equation}

A based submodule of a based  $R$-module is defined in an obvious
way. Based modules enjoy a complete reducibility \cite[Lemma
1]{O}; i.e., every based $R$-submodule of a based $R$-module $M$  
has a direct complement  which is also a based $R$-submodule of $M$.

A based ring is {\em pointed} if all its basic elements are invertible.
The invertible elements of $R$ generate a maximal pointed based subring
of $R$ which we will denote by $R_{pt}$. A typical example of a pointed
based ring is $(\mathbb{Z}G,\, G)$, where $G$ is a finite group.

Let $R$ be a based ring and let $R_+$ consist of all
$\mathbb{Z}_+$-linear combinations of  $\{X_i\}_{i\in
\mathcal{I}}$. We will say that $X\in R_+$ contains $Y\in R_+$ if
$X-Y \in R_+$. A similar terminology will be used for based
modules.

Let $G$ be a finite group. We say that a based ring $(R,\, B)$ is {\em graded} by $G$ if there is
a partition $B= \cup_{g\in G} B_g$ such that $R =\oplus_{g\in G}\, R_g$, where $R_g$ is
a $\mathbb{Z}$-submodule of $R$ generated by $B_g$ and $R_gR_h \subset R_{gh},\, R_g^* = R_{g^{-1}}$
for all $g,h\in G$.  Such a grading is {\em faithful} if $R_g\neq 0$ for all $g\in G$.

\subsection{Frobenius-Perron dimensions}

Let us recall the definition and basic properties of Frobenius-Perron dimensions
in based rings and their based modules from \cite[Section 8.1]{ENO}.  
We note that Frobenius-Perron dimensions 
for commutative based rings were defined 
and used in the book \cite{FK}.

Let $R$  be a unital based ring with basis $B= \{X_i\}_{i\in\mathcal{I}}$.
Let $X$ be a non-zero element of $R_+$ and let 
$N_X$ be the matrix of multiplication by $X$ in the basis
$B$. This matrix has nonnegative entries and is not nilpotent. 
The Frobenius-Perron theorem (see \cite{G}) implies that $N_X$ has a 
positive eigenvalue. The largest such eigenvalue is called 
the {\em Frobenius-Perron dimension} of $X$ and is denoted by $\FPdim(X)$.

The assignment $X \mapsto \FPdim(X)$ extends to  
a homomorphism from $R$ to $\mathbb{R}$. This is the unique such homomorphism 
with the property that it sends basic elements to positive numbers.
Clearly, Frobenius-Perron dimensions of elements of $R$ are algebraic integers.
The Frobenius-Perron dimension of $R$ is defined to be $\FPdim(R) =\sum_{i\in \mathcal{I}}\, \FPdim(X_i)^2$. 

Let $R =\oplus_{g\in G}\, R_g$ be a faithfully graded based ring 
with basis $B= \cup_{g\in G} B_g$.  
Define $\FPdim(R_g) =\sum_{X_i\in B_g}\, \FPdim(X_i)^2$.  
Then  $\FPdim(R_g) = \FPdim(R)/|G|$ for all $g\in G$ \cite[Proposition 8.20]{ENO}.
In particular, $|G|$ divides  $\FPdim(R)$ in the ring of algebraic integers.

Let $M$ be an indecomposable based module over a based ring $R$ with a basis $\{V_j\}_{j\in
\mathcal{J}}$. 
For every $X\in R$ let $L_X$ denote the matrix of left multiplication
by $X$ in this basis. There exists a unique (up to a positive scalar multiple)
common eigenvector $Q=~\sum_{j\in\mathcal{J}}\, \mu_j V_j \in M\otimes_\mathbb{Z} \mathbb{R}$ 
of the matrices $L_X,\, X\in R$, such that $\mu_j>0$ for all $j \in \mathcal{J}$
\cite[Proposition 8.5]{ENO}. The corresponding eigenvalue of $L_X$ is  $\FPdim(X)$.
Define the Frobenius-Perron dimension of the basic element $V_j\in M$ by $\FPdim(V_j) =\mu_j$.
Unlike the Frobenius-Perron dimensions of elements
of $R$, the dimensions  $\FPdim(V_j),\, j\in\mathcal{J}$, are defined only up to a common scalar. 
The same is true for the Frobenius-Perron dimension of $M$ defined by
$\FPdim(M) = \sum_{j\in\mathcal{J}}\, \FPdim(V_j)^2$.
We note, however, that the ratios $\FPdim(M)/\FPdim(V_j)^2$ are well defined
since they {\em do not} depend on scalings of $Q$. This fact will be used in
Section~\ref{dvisibility section}. Note also that for all $X\in R$ and $V\in M$ we have
$\FPdim(XV) = \FPdim(X)\FPdim(V)$; i.e.,  any choice of Frobenius-Perron
dimensions in $M$ defines an $R$-module homomorphism from $M$ to $\mathbb{R}$.

\subsection{Fusion categories}

By a {\em fusion category} we mean a $k$-linear semisimple rigid tensor category $\C$ with finitely
many isomorphism classes of simple objects, finite dimensional spaces of morphisms,  and such that
the unit object of $\C$ is simple.
We refer the reader to \cite{ENO} for a general theory of such categories.

Examples of fusion categories include representation categories of semisimple
Hopf and quasi-Hopf algebras as well as semisimple tensor categories 
coming from quantum groups
and affine Lie algebras \cite{BK} and from subfactor theory.

Let $\C$ be a fusion category. Its Grothendieck ring $K_0(\C)$ is the free $\mathbb{Z}$-module
generated by the isomorphism classes of simple objects of $\C$ with the multiplication
coming from the tensor product in $\C$. The Grothendieck ring of a fusion category
is a based unital ring.  The Frobenius-Perron dimensions of objects 
in $\C$ (respectively, $\FPdim(\C)$)
are defined as the Frobenius-Perron dimensions of their images in the based ring $K_0(\C)$
(respectively, as $\FPdim(K_0(\C))$). For a semisimple quasi-Hopf algebra 
$H$ one has $\FPdim(X) = \dim_k(X)$ for all
$X$ in $\Rep(H)$, and so $\FPdim(\Rep(H)) = \dim_k(H)$.

By a {\em fusion subcategory} of a fusion category $\C$ we understand a full tensor subcategory of $\C$.
An example of a fusion subcategory is the maximal {\em pointed} subcategory $\C_{pt}$ generated
by the invertible objects of $\C$.

A fusion category $\C$ is {\em pseudounitary} if its categorical dimension $\dim(\C)$ 
coincides  with its Frobenius-Perron dimension, see \cite{ENO} for details. In this case
$\C$ admits a canonical spherical structure (a tensor  isomorphism between the identity
functor of $\C$ and the second duality functor) with respect to which categorical
dimensions of objects coincide with their Frobenius-Perron dimensions
\cite[Proposition 8.23]{ENO}. We will use this fact in Section~\ref{modular stuff}.

If we identify (the isomorphism classes of ) objects $X,Y\in \C$ with the corresponding elements
of  $K_0(\C)$ then $(X,\, Y) =\dim_k\, \Hom_\C(X, Y)$.

A fusion category $\C$ is graded by a finite group $G$ if the based ring $K_0(\C)$ is $G$-graded.
That is, $\C$ decomposes into a direct sum of full abelian subcategories
$\C =\oplus_{g\in G}\,\C_g$ such that $\C_g^* = \C_{g^{-1}}$
and the tensor product maps $\C_g\times \C_h$ to $\C_{gh}$ for all $g,h\in G$.

\section{Graded  based rings and fusion categories}

\subsection{The adjoint subcategory}
Let $R$ be a based ring with a $\mathbb{Z}_+$-basis $B =\{X_i\}_{i\in \mathcal{I}}$.
Let $R_{ad}$ be the minimal based subring of $R$ with the property
that $X_i X_i^*$ belongs to $R_{ad}$ for all $i\in \mathcal{I}$; i.e.,
$R_{ad}$ is  generated by all basic elements of $R$ contained in $X_i X_i^*,\,
i\in\mathcal{I}$.

Equivalently, let us define $I(1) := \sum_{i\in \mathcal{I}}\, X_iX_i^*$.
Then $R_{ad}$ is the $\mathbb{Z}$-linear span of 
basic elements contained in $I(1)^n,\, n=1,2,\dots$. Note that
$1\in R_{ad}$.

\begin{remark}
\begin{enumerate}
\item[(1)]
Here the index ${ad}$ stands for ``adjoint'' which is justified by
the following observation. If $H$ is a semisimple Hopf algebra and
$R =K_0(\Rep(H))$ is the Grothendieck ring of the representation
category of $H$ then  $R_{ad}$  is the subring generated by the
subrepresentations of the adjoint representation of $H$.
\item[(2)]
Recall from \cite[Section 5.8]{ENO} that for a fusion category $\C$
with simple objects $\{X_i\}_{i\in \mathcal{I}}$ the induction functor $I: \C \to \Z(\C)$
to the center of $\C$ is defined as the left adjoint of the forgetful functor
$F: \Z(\C)\to \C$. One has $F(I(1)) \cong \sum_{i\in \mathcal{I}}\, X_i\ot X_i^*$, which explains
our notation.
\end{enumerate}
\end{remark}

For a fusion category  $\C$ let $\C_{ad}$ be the full tensor subcategory
of $\C$ generated by all subobjects of $X\ot X^*$ where $X$ runs
through simple objects of $\C$.  We have $K_0(\C_{ad}) = K_0(\C)_{ad}$.

\begin{example}
Let $G$ be a finite group and let $\C=\Rep(G)$ be the representation category of $G$.
Then $\C_{ad} = \Rep(G/Z(G))$, where $Z(G)$ is the center of $G$.
\end{example}

\begin{proposition}
\label{I1 is central} 
\begin{enumerate}
\item[(1)]
The element $I(1)\in  R_{ad}$ is central in $R$.
\item[(2)]
A based $R_{ad}$-subbimodule $M$ of $R$ is indecomposable 
if and only if it is indecomposable as a left (or right) $R_{ad}$-module.
\end{enumerate}
\end{proposition}
\begin{proof}
To prove (1) we first observe the equality  
\begin{equation}
\label{I1 central}
\sum _{i\in \mathcal{I}}\, X_k X_i\otimes_\mathbb{Z}  X_i^* 
= \sum _{i\in \mathcal{I}}\, X_i \otimes_\mathbb{Z} X_i^* X_k, \qquad \mbox{ for all } k\in\mathcal{I},
\end{equation} 
which can be obtained by evaluating the pairing $(,)$
of both sides with \linebreak $X_j \otimes_\mathbb{Z} X_l,\, j,l\in~\mathcal{I}$, and using 
Definition~\ref{based ring}. Applying the multiplication of $R$ to both sides of
\eqref{I1 central} we obtain $X_k I(1) = I(1)X_k$ for all $k\in\mathcal{I}$.

Since every basic element of $R_{ad}$ is contained in some power of $I(1)$, 
every one-sided based $R_{ad}$-submodule of $M$ is automatically an $R_{ad}$-subbimodule.
Indeed, if $XY \in M$ for all $X\in R_{ad}$ and $Y\in M$ then $I(1)^n Y = Y I^n(1)\in M$
for all non-negative integers $n$, and hence $YX\in M$. 
\end{proof}

\subsection{The universal grading of a based ring and
of a fusion category}
Let $R$ be a based ring.
We can view $R$ as a based  $R_{ad}$-bimodule. As such, it decomposes
into a direct sum of indecomposable based $R_{ad}$-bimodules:
$R = \oplus_{a\in A}\, R_a$,
where $A$ is the index set.
This decomposition is unique up to a permutation of $A$.
We may assume that there is an element $1\in A$ such that $R_1 = R_{ad}$.
Note that $(R_a)^*=\{X^* \mid X \in R_a\},\, a\in A,$
is an indecomposable based $R_{ad}$-submodule of $R$ and hence $(R_a)^* = R_{a^*}$ for some
$a^*\in A$.

\begin{lemma}
\label{inverse}
For all $X_a, Y_a \in R_a,\, a\in A$, we have $X_aY_a^*\in R_{ad}$.
\end{lemma}
\begin{proof}
Recall that each $X\in R_{ad}$ is contained in some power of $I(1)$.
Observe that $M= \{X\in R_a \mid X \mbox{ is contained in } I(1)^nY_a \mbox{ for some } n\}$
is a $R_{ad}$-submodule of $R_a$. 
Since  $R_a$ is an indecomposable based $R_{ad}$-bimodule,
it follows from Proposition~\ref{I1 is central}(2) that $M=R_a$ and therefore
$X_a$ is contained in $I(1)^nY_a$ for some $n$.
Hence, $X_aY_a^*$
is contained in $I(1)^n Y_a Y_a^* \in R_{ad}$, as required.
\end{proof}

\begin{theorem}
\label{constructing a group}
There is a canonical group structure on the index set $A$
with the multiplication defined by the following property:
\begin{equation}
\label{product in A}
ab =c
 \quad \mbox{if and only if } X_a X_b \in R_{c} ,\quad \mbox{for all }
X_a\in R_a, X_b\in R_b,\, a,b,c\in A.
\end{equation}
The identity of $A$ is $1$ and the inverse of $a\in A$ is $a^*$.
\end{theorem}
\begin{proof}
We need to check that the binary operation in \eqref{product in A} is well defined.
Let $a,b\in A$ and let $X_a, Y_a \in R_a, X_b, Y_b\in R_b$ be basic elements of $R$.

Suppose that the product
$X_aX_b$ contains a basic element $X_c\in R_c$
and the product $Y_aY_b$ contains a basic element $Y_d\in R_d$
for $ c\neq d$.
By Lemma~\ref{inverse},
there is a positive integer $n$ such that   the element
$X_a I(1)^n Y_a^*\in R_1$ contains $X_a X_b Y_b^* Y_a^*   = (X_aX_b)(Y_aY_b)^*$ and, therefore,
$Y= X_c Y_d^*$ is in $R_1$. Multiplying both sides of the last equality  by $Y_d$
on the right  we conclude that $YY_d$ is in $R_c\cap R_d$, a contradiction.

Thus, $X_aX_b$ and $Y_aY_b$ both belong to the same component $R_c$
and so  the binary operation \eqref{product in A} is well-defined.
It is easy to see that it defines a group structure on~$A$.
\end{proof}

\begin{definition}
We will call the grading $R =\oplus_{a\in A}\, R_a$ constructed in Theorem~\ref{constructing a group}
the {\em universal grading} of $R$.  The group $A$ will be called the {\em universal grading group} of $R$
and denoted by $U(R)$.
\end{definition}

\begin{corollary}\label{faithful}
\label{universality}
Every based ring $R$ has a canonical faithful grading by the group $U(R)$. Any other faithful grading
of $R$ by a group $G$ is determined by a surjective group homomorphism $\pi : U(R)\to G$.
\end{corollary}
\begin{proof}
Let $R =\oplus_{g\in G}\, R^g$ be a grading of $R$.  Since for every basic $X\in R$ we have $XX^*\in R^1$
it follows that  $R^1$ contains $R_{ad}$ as a based subring. Hence, each $R^g$ is a based  $R_{ad}$-submodule
of $R$. This means that every component  $R_a,\, a\in U(R)$, of the universal grading $R =\oplus_{a\in U(R)}\, R_a$
of $R$ belongs to some $R^{\pi(a)}$ for some well-defined $\pi(a)\in G$.  Clearly, the map $a \mapsto \pi(a)$
is a surjective homomorphism.
\end{proof}

For a fusion category $\C$ let $U(\C)=U(K_0(\C))$. We will call $U(\C)$ the {\em universal
grading group} of $\C$.

\begin{theorem}\label{hopfad}
Let $H$ be a semisimple Hopf algebra and let $\C=\Rep(H)$.  There exists
a unique Hopf subalgebra $K$ of $H$ which is contained in the center of $H$
and is maximal with respect to this property. Let $H_{ad} := H/HK^+$
be the quotient of $H$ by $K$ (see \cite{Mo}). Then $\C_{ad} = \Rep(H_{ad})$.
Furthermore, $K = \Fun(U(\C))$. 
\end{theorem}
\begin{proof}
Clearly such $K$ exists. 
If $K_1, K_2$ are Hopf subalgebras of $H$ contained in its center, then so is
the Hopf subalgebra $K_1 K_2$, so that $K$ is unique. 
Let $\Delta$ and $S$ denote the comultiplication and antipode of $H$.
Every grading of $\Rep(H)$ 
comes from a decomposition $H = \oplus_{a \in A}\, H_a$ where each $H_a$ is an
ideal of $H$ such that  $\Delta(H_a) \subseteq \bigoplus_{xy =a}\, H_x \ot H_y$ 
and $S(H_a) = H_{a^{-1}}$. Let $p_a\in H$ be the central idempotents such that
$p_a H = H_a,\, a\in A$. Then $\Delta(p_a) = \sum_{xy =a}\, p_x \ot p_y$; i.e., 
the linear span of $p_a, \, a\in A$, is 
a Hopf subalgebra of $H$ (isomorphic to $\Fun(A)$)
and is contained in the center of $H$. Clearly, the universal grading of $\Rep(H)$
corresponds to a maximal such subalgebra. The trivial component of $H$ is
precisely the quotient $H/HK^+$, so that $ \C_{ad} = \Rep(H_1) = \Rep(H/HK^+)$.
\end{proof}

Let $\C$ be a fusion category and
let $\Aut_\otimes(\id_\C)$ denote the group of tensor automorphisms
of the identity functor of $\C$. Any $\Phi \in \Aut_\otimes(\id_\C)$
is determined by a collection of non-zero scalars $\{\Phi(X)\}_X$,
one for every simple object $X$ in $\C$, such that
$\Phi(X_1)\Phi(X_2) = \Phi(Y)$ whenever $Y$ is contained
in $X_1\otimes X_2$. Clearly, $\Aut_\otimes(\id_\C)$ 
is an abelian group.

For every abelian group $A$ let $\widehat{A}$ denote the group of characters of $A$.

\begin{proposition}\label{maxab}
\label{abelian case}
Let $\C$ be a fusion category and let $G=U(\C)$ be its universal grading group.
Let $G_{ab}$ be the maximal abelian quotient group of $G$. There is an isomorphism
between $\widehat{G_{ab}}$ and the group $\Aut_\otimes(\id_\C)$.
\end{proposition}
\begin{proof}
Define a grading of $\C$ by $\widehat{\Aut_\otimes(\id_\C)}$
by setting
\begin{equation}
\label{grading by Aut}
\C_\chi= \{X \in \C \mid \Phi_X =\chi(\Phi)\id_X\mbox { for all } \Phi \in \Aut_\otimes(\id_\C) \},
\end{equation}
where $\chi \in \widehat{\Aut_\otimes(\id_\C)}$. Since $G_{ab}$ is the
universal {\em abelian} grading group of $\C$ it follows that
there is a surjective group homomorphism $G_{ab} \to \widehat{\Aut_\otimes(\id_\C)}$.

Conversely, let $\C =\oplus_{a\in G_{ab}}\, \C_a$ be the grading of $\C$ by $G_{ab}$
constructed in Corollary~\ref{universality}.
For each $\chi \in \widehat {G_{ab}}$ define $\Phi_\chi \in \Aut_\otimes(\id_\C)$
by  $(\Phi_\chi)_X = \chi(a) \id_X$ for all $X\in \C_a$. The map
$\chi\mapsto \Phi_\chi$ is an injective group homomorphism, and hence it establishes
the required isomorphism.
\end{proof}

\subsection{Based rings of integer Frobenius-Perron dimension}

It was shown in \cite[Section 8.5]{ENO} that a fusion category of integer
Frobenius-Perron dimension is  pseudounitary (and hence, admits a pivotal structure 
with respect to which
the categorical dimensions of all simple objects coincide with their Frobenius-Perron 
dimensions).

Let us analyze the structure of  based rings of integer Frobenius-Perron dimension.

\begin{theorem}\label{int}
\label{when FPdim in Z}
Let $R$ be a based ring such that $\FPdim(R) \in \mathbb{Z}$. Then there is an elementary abelian
$2$-group $E$, a set of distinct square free positive integers $n_x,\, x\in~E$, with $n_0=1$,
and a faithful grading $R =\oplus_{x\in E}\, R(n_x)$
such that  $\FPdim(X) \in \mathbb{Z}\sqrt{n_x}$ for each $X\in R(n_x)$.
\end{theorem}
\begin{proof}
By \cite[Proposition 8.27]{ENO} every basic element of $R$ has dimension $\sqrt{N}$ for some $N\in \mathbb{Z}$.
Let $R(1) \subset R$ be the based subring of $R$ generated by all basic elements of integer dimension.
Observe that $R_{ad} \subset R(1)$ and for each square free $n\in \mathbb{Z}$ the basic elements of $R$
whose dimension is in $\mathbb{Z}\sqrt{n}$ generate an $R(1)$-subbimodule $R(n)$ of $R$.  Let
\[
E = \{n  \mbox{ is square free } \mid  \exists X\in R_+,\, X\neq 0, \mbox{ such that } \FPdim(X)\in
\mathbb{Z}\sqrt{n}\}.
\]
It is clear that for $X\in R(n)$
and $Y\in R(m)$ their product $XY$ is in $R((nm)')$ where $l'$ denotes the square free part of $l$.
This defines a commutative group operation on $E$ and a grading on $R$. Since the order of every $e\in E$ is
at most two, $E$ is an elementary abelian $2$-group.
\end{proof}

\begin{corollary}\label{divprop}
Let $\C$ be a fusion category of odd Frobenius-Perron dimension. Then the 
Frobenius-Perron dimension of every object of $\C$ is an integer 
(and hence, $\C$ is the category of representations of some
semisimple quasi-Hopf algebra, see \cite{ENO}).
\end{corollary}


\section{The central series and nilpotency for based rings
and fusion categories}

\subsection{The upper central series}
Let $R$ be a based ring. Let $R^{(0)}= R,\, R^{(1)} = R_{ad},$
and $R^{(n)} = (R^{(n-1)})_{ad}$ for every integer $n\geq 1$.

\begin{definition}
\label{upper central series}
The non-increasing sequence of based subrings of $R$
\begin{equation}
R= R^{(0)} \supseteq  R^{(1)} \supseteq \cdots \supseteq R^{(n)}
\supseteq \cdots
\end{equation}
will be called the {\em upper central series} of $R$.
\end{definition}

Similarly, for a fusion category $\C$ we define
$\C^{(0)}= \C,\, \C^{(1)} = \C_{ad},$
and $\C^{(n)} = (\C^{(n-1)})_{ad}$ for every integer $n\geq 1$.

\begin{definition}\label{upper central series 1}
\label{upp central series}
The non-increasing sequence of fusion subcategories of $\C$
\begin{equation}
\C= \C^{(0)} \supseteq  \C^{(1)} \supseteq \cdots
\supseteq \C^{(n)} \supseteq \cdots
\end{equation}
will be called the {\em upper central series} of $\C$.
\end{definition}

\begin{example}
For every group $H$ let $Z(H)$ denote its center.
Let $G$ be a finite group and  $\C=\Rep(G)$. Let
\[
\{1\} = C^0(G) \subseteq C^{1}(G) \subseteq
\cdots  \subseteq   C^{n}(G) \subseteq \cdots
\]
be the upper central series of $G$; i.e.,
$C^0(G):=\{1\}, C^{1}(G):= Z(G)$ and for $n\geq 1$
the subgroup  $C^n(G)$
is defined by $C^n(G)/C^{n-1}(G) = Z(G/C^{n-1}(G))$.
Then $\C^{(n)} = \Rep(G/C^n(G))$,
so that our definition of the upper central series
agrees with the classical one.
\end{example}

\begin{definition}
\label{nilpotent}
A based ring $R$ is {\em nilpotent} if its upper central series
converges to $\mathbb{Z}1$; i.e., $R^{(n)} = \mathbb{Z}1$
for some $n$. The smallest number $n$ for which this happens
is called the {\em nilpotency class} of $R$.

A fusion category $\C$ is {\em nilpotent} if
its upper central series converges to $\Vec$; i.e.,
$\C^{(n)} =\Vec$ for some $n$. The smallest such $n$
is called the  {\em nilpotency class} of $\C$.

A semisimple Hopf algebra (or  a quasi-Hopf algebra, or
a weak Hopf algebra) is {\em nilpotent}
if its representation category is nilpotent.
\end{definition}

\begin{example}
\begin{enumerate}
\item[(1)]
By \cite[Theorem 8.28]{ENO}, every fusion category whose dimension is a prime power
(i.e., a {\em $p$-fusion category}) is graded by $\mathbb{Z}/p\mathbb{Z}$, and hence
is nilpotent. In particular,
every semisimple Hopf or quasi-Hopf algebra of a prime power dimension
is nilpotent.
\item[(2)]
Semisimple Hopf algebras of dimension $2n^2,\, n\geq 2$,
constructed by G.~Kac and V.~Paljutkin in \cite{KP} (see also \cite{Se})
are nilpotent.
These Hopf algebras have algebra structure $k^{n^2} \oplus M_n(k)$
and representation categories of Tambara-Yamagami type \cite{TY}.
In particular, a semisimple nilpotent Hopf algebra
is not necessarily a tensor product of Hopf algebras of prime power dimension,
which is in contrast with a classical result in group theory.
\end{enumerate}
\end{example}

Let $\C$ and $\D$ be fusion categories.
Recall that for a tensor functor $F: \C\to \D$ its image 
$F(\C)$ is the fusion subcategory of $\D$ generated by all simple objects
$Y$ in $\D$ such that $Y \subseteq F(X)$ for some simple $X$ in $\C$.
The functor $F$ is called {\em surjective} if $F(\C) =\D$.

\begin{proposition}
\label{subcat+ surj image}
Let $\C$ be a nilpotent fusion category.
Every fusion subcategory $\mathcal{E}\subset \C$
is nilpotent. If $F: \C\to \D$
is a surjective tensor functor, then $\D$ is nilpotent.
\end{proposition}
\begin{proof}
We have $\mathcal{E}^{(n)} \subset  \mathcal{C}^{(n)}$ and
$\mathcal{D}^{(n)}\subset F(\mathcal{C}^{(n)})$ for all $n$.
\end{proof}

\begin{remark}
\label{long remark}
\begin{enumerate}
\item[(1)]
Let $G$ be a finite group and  $\C=\Rep(G)$. Then
$\C$ is nilpotent if and only if $G$ is nilpotent.
\item[(2)]
By Theorem~\ref{constructing a group},
a based ring $R$ (respectively,
a fusion category $\C$)
is nilpotent if and only if every non-trivial
based subring of $R$ ( respectively,
a non-trivial fusion subcategory of $\C$)
has a non-trivial grading.
\item[(3)]
A fusion category $\C$ is nilpotent if and only
if $K_0(\C)$ is nilpotent. The nilpotency
class of $\C$ is equal to the nilpotency class of $K_0(\C)$.
\item[(4)]
Pointed based rings (respectively, pointed fusion categories)
are precisely the nilpotent
based rings (respectively, fusion categories) of nilpotency class~$1$.
Fusion categories constructed by  Tambara and Yamagami
\cite{TY} have nilpotency class $2$.
\item[(5)]
A nilpotent semisimple Hopf algebra is lower solvable
in the sense of \cite{MW}. (The special case of a twisted 
group algebra of a nilpotent group is discussed in \cite{GN}.)
\item[(6)]
It follows from Theorem~\ref{constructing a group}
that a nilpotent fusion
category comes from  a sequence of gradings, in particular
it has an integer Frobenius-Perron dimension.
It follows from results of \cite{ENO}
that a nilpotent fusion category
$\C$ is pseudounitary,
i.e., its global dimension coincides with $\FPdim(\C)$,
and that $\C$ admits a pivotal (in fact, spherical) structure.
\end{enumerate}
\end{remark}

\subsection{Lower central series  for commutative based rings}
For a commutative based ring $R$ we can define
the notion of a lower central series.

\begin{definition}\label{com}
Let $S$ be a based subring of a commutative based ring $R$.
We define the {\em commutator} $S^{co}$ of
$S$ in $R$ to be the  based subring of
$R$ generated by all basic elements $Y\in R$ such that
$YY^*\in S$.
\end{definition}

\begin{remark}
Note that the $\mathbb{Z}$-span of basic elements $Y\in R$
with the property $YY^*\in S$ is already a based subring of $R$
and so is equal to $S^{co}$.
Indeed, if $Y_1Y_1^*, Y_2Y_2^* \in S$ and $Y$ is contained
in $Y_1Y_2$ then $YY^*$ is contained in $(Y_1Y_1^*)(Y_2Y_2^*) \in~S$.
\end{remark}

\begin{definition}\label{comfus}
Let $\C$ be a fusion category such that $K_0(\C)$ is commutative
(e.g., $\C$ is braided).
Let $\mathcal{K}$ be a fusion subcategory  of $\C$.
Define the {\em commutator} $\mathcal{K}^{co}$ of $\mathcal{K}$
to be the  fusion subcategory
of $\C$ generated by all simple objects $Y$ of $\C$ such that $Y\otimes Y^*\in \mathcal{K}$.
\end{definition}
One has  $K_0(\mathcal{K}^{co}) = K_0(\mathcal{K})^{co}$.

\begin{example}
\label{about commutators}
Let $G$ be a finite group and  let $\C=\Rep(G)$.
Any fusion subcategory  of $\C$ is of the form
$\Rep(G/N)$ for some normal subgroup $N$ of $G$.
The simple objects of the
category $\Rep(G/N)^{co}$  are irreducible 
representations $Y$ of $G$ for which $ Y \otimes Y^*$  restricts
to the trivial representation of $N$.
Equivalently, each $x\in N$ acts on $Y$ by  a scalar.
This happens if and only if the commutator group
$[G, N] := \langle gxg^{-1}x^{-1}\mid g\in G, x \in N \rangle $
acts trivially on $Y$.  Thus, $\Rep(G/N)^{co} = \Rep(G/[G, N])$.
\end{example}

For a commutative based ring $R$ we define
$R_{(0)}= \mathbb{Z}1$,  $R_{(1)} = R_{pt}$
(the based subring of $R$ spanned by the invertible basic objects of $R$),
and $R_{(n)} = (R_{(n-1)})^{co}$ for every integer $n\geq 1$.

\begin{definition}
\label{lower central series}
Let $R$ be a commutative based ring.
The non-decreasing sequence of based subrings of $R$
\begin{equation}
\mathbb{Z}1 = R_{(0)} \subseteq  R_{(1)} \subseteq \cdots
\subseteq R_{(n)} \subseteq \cdots
\end{equation}
will be called the {\em lower central series} of $R$.
\end{definition}

Similarly, for a fusion category $\C$ we define
$\C_{(0)}= \Vec,\, \C_{(1)} = \C_{pt}$
(the maximal pointed subcategory of $\C$)
and $\C_{(n)} = (\C_{(n-1)})^{co}$ for every integer $n\geq 1$.

\begin{definition}
\label{low central series}
The non-decreasing sequence of fusion subcategories of $\C$
\begin{equation}
\Vec= \C_{(0)} \subseteq  \C_{(1)} \subseteq \cdots \subseteq
\C_{(n)} \subseteq \cdots
\end{equation}
will be called the {\em lower central series} of $\C$.
\end{definition}

\begin{example}
Let $G$ be a finite group and  $\C=\Rep(G)$. Let
\[
G = C_0(G) \supseteq C_{1}(G) \supseteq   \cdots  \supseteq
C_{n}(G) \supseteq \cdots
\]
be the lower central series of $G$; i.e.,
$C_n(G) = [G, C_{n-1}(G)]$ for all $n\geq 1$.
Then $\C_{(n)} = \Rep(G/C_n(G))$, so that our definition
of the lower central series agrees with the classical one.
\end{example}

For a commutative based ring $R$
we can prove  an analogue of the classical result in group theory
saying that the upper central series of a group $G$
converges to $G$ if and only if
its lower central series converges to $\{ 1\}$, which gives
an alternative criterion of nilpotency \cite{H}.

\begin{lemma}
\label{coad in adco}
Let $R$ be a commutative based ring and let $S\subseteq R $ be a based subring.
Then $(S^{co})_{ad} \subseteq S \subseteq (S_{ad})^{co}$.
\end{lemma}
\begin{proof}
Both inclusions are immediate from the definitions of
the adjoint subring and  commutator.
\end{proof}

\begin{theorem}\label{upp=low}
Let $R$ be a commutative based ring and let
\begin{gather}
R = R^{(0)} \supseteq  R^{(1)} \supseteq \cdots
\supseteq R^{(n)} \supseteq \cdots \\
\mathbb{Z}1 = R_{(0)} \subseteq  R_{(1)} \subseteq \cdots
\subseteq R_{(n)} \subseteq \cdots
\end{gather}
be the upper and lower central series of $R$.
Then $R^{(n)} = \mathbb{Z}1$ if and only if
$R_{(n)} = R$, and when this is the case one has $R^{(k)}
\subseteq R_{(n-k)}$ for all $k = 0,\dots n$.
\end{theorem}
\begin{proof}
Suppose $R_{(n)}=R$; i.e., $R^{(0)} = R_{(n)}$.
If  $R^{(k)} \subseteq R_{(n-k)}$ for some $k$
($0\leq k < n$) then, using Lemma~\ref{coad in adco},
we obtain $R^{(k+1)} = (R^{(k)})_{ad} \subseteq (R_{(n-k)})_{ad}
= ({R_{(n-k-1)}}^{co})_{ad} \subseteq R_{(n-k-1)}$. Hence,
$R^{(n)} = \mathbb{Z}1$ by induction.

Conversely, suppose $R^{(n)} = \mathbb{Z}1$; i.e.,
$R^{(n)} = R_{(0)}$. If  $R^{(n-k)} \subseteq R_{k}$ for some $k$
($n\geq k > 0$) then
$R^{(n-k-1)} \subseteq ({R^{(n-k-1)}}_{ad})^{co} =
(R^{(n-k)})^{co} \subseteq (R_{(k)})^{co} = R_{(k+1)}$.
Hence, $R_{(n)}=R$.

This means that neither central series can be longer than
the other, as required.
\end{proof}

\section{Divisibility properties of Frobenius-Perron dimensions of basic
elements in a nilpotent based ring}
\label{dvisibility section}

Let $(R,\, \{X_i\}_{i\in \mathcal{I}})$ be a based ring  with the
universal grading $R =\oplus_{a\in U(R)}\, R_a$, 
and let $(M,\, \{V_j\}_{j\in \mathcal{J}})$ be an
indecomposable based $R$-module. Then $M$ is also a left based
$R_{ad}$-module and therefore decomposes into a direct sum of
indecomposable based $R_{ad}$-submodules:  $M =\oplus_{x\in S}\,
M_x$. This decomposition is unique up to a permutation of $S$.

Let us fix a Frobenius-Perron eigenvector 
$Q= \sum_{j\in \mathcal{J}}\,  \FPdim(V_j) V_j \in M \otimes_\mathbb{Z} \mathbb{R}$.
For each $x\in S$ let $\FPdim(M_x) = \sum_{V_j \in M_x}\, \FPdim(V_j)^2$.
The vector $Q$ is unique up to a scalar, hence the ratios
$\FPdim(M_x)/\FPdim(M_y)$ are well-defined for all $x,y\in S$.

\begin{proposition}
\label{constructing a Aset}
There is a canonical structure of a transitive $U(R)$-set on the index set $S$
defined by the following property:
\begin{equation}
\label{action of A}
ax =y
 \quad \mbox{if and only if } X_a V_x \in M_y ,\quad \mbox{for all }
X_a\in R_a, V_x\in M_x,\, a\in U(R),\, x,y\in S.
\end{equation}
Furthermore, $\FPdim(M_x)/\FPdim(M_y) =1$ for all $x,y\in S$.
\end{proposition}
\begin{proof}
We need to check that the action in \eqref{action of A} is well
defined. Let $X_a, Y_a\in R_a,\,a\in U(R)$, and let $V_x, U_x$ be
basic elements in $M_x,\, x\in S$. Suppose that the product
$X_aV_x$ contains a non-zero element $V_y\in M_y$ and the product
$Y_aU_x$ contains a non-zero element $U_z\in M_z$ for some $y\neq z$ in
$S$. Let $I(1) = \sum_{i\in \mathcal{I}}\, X_iX_i^* \in R_{ad}$.
Since $M_x$ is an indecomposable
based $R_{ad}$-module, for some $n\in \mathbb{Z}_+$ the element
$U_x$ is contained in $I(1)^nV_x$. Therefore, we have
\[
(V_x,\, X_a^*V_y) = (X_aV_x,\, V_y)>0
\]
 and
\[
(V_x,\, I(1)^{n} Y_a^*U_z) = (Y_aI(1)^nV_x,\, U_z)\geq (Y_aU_x, U_z) >0.
\]
Combining these inequalities we get
$(Y_aX_a^*V_y,\, I(1)^n U_z) > 0$. Since $Y_a X_a^* \in R_{ad}$,
we obtain a contradiction: $(V'_y,\, U'_z)>0$ for $V'_y= Y_aX_a^*V_y \in M_y$
and $U'_z=  I(1)^n U_z \in M_z$, while $y\neq z$.  This proves that the action in question
is well-defined.

To prove the statement about Frobenius-Perron  dimensions let
$$E= \sum_{a\in U(R)}\, E_a \in R\otimes_\mathbb{Z} \mathbb{R}$$ be the uniquely determined
virtual regular element of $R$, where  $$E_a =\sum_{X_i\in R_a}\, \FPdim(X_i) X_i,$$
and let $Q= \sum_{x\in S}\, Q_x\in M \otimes_\mathbb{Z} \mathbb{R}$, where  
$$
Q_x   =\sum_{V_j\in M_x}\,\FPdim(V_j)V_j  \in M_x \otimes_\mathbb{Z} \mathbb{R}
$$
for each $x\in S$.

The positive number $\FPdim(E_a), \,a\in U(R),$  
does not depend on $a$ by \cite[Proposition 8.20]{ENO}.
Let us denote it by $d$. One has $E_aQ = d Q$, whence $E_aQ_x = dQ_{ax}$. 
Since $\FPdim: M \to \mathbb{R}$ is an $R$-module homomorphism, taking $\FPdim$
of both sides of the last equality we obtain
$\FPdim(Q_x) = \FPdim(Q_{ax})$.  Hence,  $ \FPdim(M_x) = \FPdim(M_{ax}).$ 
Since the action of $U(R)$ on $S$ is transitive, the statement follows.
\end{proof}

\begin{theorem}
\label{inductive step} 
Let $R$ be a nilpotent based ring, let $M$ be an
indecomposable based $R-$module, and let $V\in M$ be a basic element.
Then $\FPdim(M)/\FPdim(V)^2 \in \mathbb{Z}$. 
In particular,
$\FPdim(R_{ad})/\FPdim(X)^2 \in \mathbb{Z}$ for every basic $X$ in
$R$.
\end{theorem}
\begin{proof}
Note that the ratio $\FPdim(M)/\FPdim(V)^2$ is well defined.

We use induction on the nilpotency class of $R$.
For pointed based rings (i.e., those of nilpotency class $1$)
the statement is clear, since any indecomposable based module $M$  over $\mathbb{Z}G$,
where $G$ is a finite group, comes from a transitive  $G$-set and one can choose
$\FPdim(V) = 1$ for any basic $V\in M$. 

Let $R$ be nilpotent of class $n> 1$,
then $R_{ad}$ is nilpotent of class $n-1$. Suppose the statement is true
for all indecomposable $R_{ad}$-modules.  Let $M$ be an indecomposable $R$-module
and let  $M_x,\, x\in S$, be its indecomposable $R_{ad}$-submodules. 

Let $V$ be a basic object of $M$, then $V\in M_x$ for some $x\in S$.
Applying Proposition~\ref{constructing a Aset}, 
and using the inductive assumption,
we have
\begin{equation*}
\frac{\FPdim(M)}{\FPdim(V)^2} = \frac{|S|\FPdim(M_x)}{\FPdim(V)^2} \in \mathbb{Z},
\end{equation*}
as required.
The second statement follows since each component $R_a$  in the universal grading 
$R~=~\oplus_{a\in U(R)}\,R_a$ 
is an indecomposable based $R_{ad}$-module and $\FPdim(R_a) =\FPdim(R_{ad})$
for all $a\in A$.
\end{proof}

\begin{corollary}
\label{nilpotent divisibility}
Let $\C$ be a nilpotent fusion category. Then for each simple object $X$ in $\C$,
$\FPdim(X)^2$ divides $\FPdim(\C_{ad})$.
\end{corollary}

\begin{remark}
\begin{enumerate}
\item[(1)]
In the case when $\C$ is the representation category of a nilpotent group, we recover a well-known
fact about group representations: the square of the degree of an irreducible representation
of a finite nilpotent  group divides the index of its center.
\item[(2)]
It follows from \cite{MW} that for an irreducible representation $V$ of an upper or lower
semi-solvable Hopf algebra $H$ one has $\dim_k(V) | \dim_k(H)$.
Corollary~\ref{nilpotent divisibility} generalizes this result in 
the special case when $H$ is nilpotent (e.g., has dimension $p^n$, $p$ a prime).

\end{enumerate}
\end{remark}

\section{The lower and upper central series of modular categories}
\label{modular stuff}

Let $\C$ be a modular category \cite{BK, T} 
with simple objects $X_i,\, i\in
\mathcal{I}$, the braiding  $c_{X,Y} : X \ot Y \cong Y\otimes X$ for all $X,Y\in \C,$
the $S$-matrix $(s_{ij})_{i,j\in
\mathcal{I}}$, where $s_{ij} =\tr_q(c_{X_i, X_j} \circ c_{X_j, X_i})$, and the
$T$-matrix $(\delta_{ij}\theta_i)_{i,j\in \mathcal{I}}$. 
Each dimension $d_i:= d(X_i) =s_{0i}$ is a  real number and each $\theta_i, \, i\in \mathcal{I}$,
is a root of unity. The
entries of the $S$-matrix are known to satisfy $s_{ij} =s_{ji}=
s_{i^*j^*}$ and $\sum_j\, s_{ij}s_{jl} = \delta_{i^*l} D$, 
where $i^*\in\mathcal{I}$ is an index such that $X_{i^*} =X_i^*$
and $D =\sum_{i\in\mathcal{I}}\, d_i^2$ is the categorical dimension of~$\C$.

A theorem  due to J.~de~Boere, J.~Goeree, A.~Coste and T.~Gannon
states that the entries of the S-matrix of a semisimple modular
category lie in a cyclotomic field, see \cite{dBG, CG}. That is,
there exists a root of unity $\xi$ such that $s_{ij}\in {\mathbb Q}(\xi)\subset \mathbb{C}$.
So we will always assume that the entries of the $S$-matrix are complex numbers.

It is known that $S/\sqrt{D}$ is a unitary matrix; i.e.,  
$s_{ij^*} =\overline{s}_{ij}$ for all $i,j\in \mathcal{I}$
\cite[Proposition 2.12]{ENO}. It follows that 
$D= \sum_i\,s_{ij}s_{i^*j}  =  \sum_i\, |s_{ij}|^2$  for each $j\in\mathcal{I}$.

Let $X_i \otimes X_j  \cong \oplus_{k \in \mathcal{I}}\, N_{ij}^k
X_k,\, i,j\in \mathcal{I}$, be the fusion rules of $\C$.
The Verlinde formula (see e.g., \cite[3.1]{BK})
relates the fusion coefficients and entries of the $S$-matrix:
\begin{equation}
\label{verlinde1}
N_{il}^j = \sum_{k\in\mathcal{I}}\,
\frac{s_{ik}s_{lk}s_{j^*k}}{d_k}.
\end{equation}
It follows from this formula that  every homomorphism from $K_0(\C)$
to $\mathbb{C}$ has the form
\begin{equation}
\label{verlinde}
h_j: X_i \mapsto \frac{s_{ij}}{s_{0j}},\quad j\in
\mathcal{I}.
\end{equation}

The entries of the $S$- and $T$-matrices
are related by the following formula
\cite[3.1.2]{BK}:
\begin{equation}
\label{sij}
s_{ij} = \theta_i^{-1}\theta_j^{-1} \sum_{k\in
\mathcal{I}}\, N_{i^*j}^k \theta_k d_k.
\end{equation}

\subsection{The universal grading group of a modular category}

\begin{lemma}
\label{invertible} Let $j\in \mathcal{I}$.  We have $|s_{ij}|
= |d_id_j|$ for all $i\in \mathcal{I}$ if and only if $X_j$ is
invertible.
\end{lemma}
\begin{proof}
Suppose $X_j$ is invertible. Then $d_j =\pm 1$ and for each $i\in
\mathcal{I}$ there is a unique $k\in \mathcal{I}$ such that
$N_{i^*j}^k\neq 0$ (in fact, $N_{i^*j}^k=1$).  Hence, $s_{ij} = \theta_i^{-1} \theta_j^{-1}
\theta_k d_k$. Since every $\theta_i,\, i\in \mathcal{I}$, is a
root of unity and $d_k = d_id_j$ we have $|s_{ij}| = |d_id_j|$.

Conversely, if $j\in \mathcal{I}$ is such that $|s_{ij}| =
|d_id_j|$ for all $i\in \mathcal{I}$ then 
$D= \sum_i\, |s_{ij}|^2 = d_j^2 D$ and hence $d_j =\pm 1$.
Therefore,
$s_{ij}= \xi_{ij} d_i$ for some $\xi_{ij}$ such that $|\xi_{ij}| =1$.
But then it follows from \eqref{verlinde} that the image of $X_j$
under any ring homomorphism $K_0(\C) \to \mathbb{C}$ has absolute value $1$.
Hence, $\FPdim(X_j) =1$ and $X_j$ is invertible.
\end{proof}

Let $U(\C)$ be the universal grading group of $\C$ and let $G(\C)$
be the group of invertible objects of $\C$. Let $X(\C)$ denote the
set of ring homomorphisms $K_0(\C)\to \mathbb{C}$.  Formula
\eqref{verlinde} gives a bijection between $X(\C)$ and the set
$\mathcal{I}$ of isomorphism classes of simple objects of $\C$.

\begin{lemma}
\label{action of Aut} The group $\Aut_\otimes(\id_\C)$ acts on
$X(\C)$ by $\Phi h (X):= \Phi(X) h(X)$ for all $\Phi \in
Aut_\otimes(\id_\C),\, h \in X(\C)$, and simple $X\in \C$.
\end{lemma}
\begin{proof}
This is a direct consequence of the definition of a tensor automorphism.
\end{proof}

\begin{theorem}
\label{Schneider's iso}
There is a canonical isomorphism $U(\C) \cong \widehat{G(\C)}$.
\end{theorem}
\begin{proof}
Since the group $U(\C)$ is abelian, it suffices to establish
an isomorphism $G(\C) \cong \Aut_\otimes(\id_\C)$,
thanks to Proposition~\ref{abelian case}.

For any $\Phi \in Aut_\otimes(\id_\C)$
the function $X_i \mapsto \Phi(X) d_i$ extends to a
homomorphism $K_0(\C)\to\mathbb{C}$
and hence there is a unique $j_\Phi\in \mathcal{I}$ such that
$h_{j_\Phi}(X_i) = \Phi(X_i) d_i$.  That is,
\begin{equation}
\label{Phi to h}
s_{ij_\Phi} =  \Phi(X_i) d_i d_{j_\Phi}, \quad
\mbox{ for all } i\in \mathcal{I}.
\end{equation}
Since $\Phi(X)$ is a root of unity we have
$|s_{ij_\Phi}| = |d_i d_{j_\Phi}|$.
By Lemma~\ref{invertible}, $X_{j_\Phi}$ is an invertible object of $\C$.

Conversely, for any invertible $X_j\in \C$ let $\Phi_j(X_i) =
\frac{s_{ij}}{d_id_j},\, i\in \mathcal{I}$. To show that $\Phi_j$
defines a tensor automorphism of $\id_\C$ we need to check that
\begin{equation}
\label{ss=s}
\frac{s_{kj}}{d_kd_j} = \frac{s_{ij}}{d_id_j} \frac{s_{lj}}{d_ld_j}
\end{equation}
for all $k, i, l\in \mathcal{I} $ such that $X_k$ is contained in
$X_i \otimes X_l$. Let $\gamma_{i},\, i\in\mathcal{I},$ be the scalar
such that $c_{X_j, X_i} \circ c_{X_i, X_j} = \gamma_{i} \id_{X_i \ot X_j}$. 
Clearly, $\gamma_i = \frac{s_{ij}}{d_id_j}$. Using the hexagon identity, we compute
\[
c_{X_i\ot X_l, X_j} \circ c_{X_j, X_i\ot X_l} = c_{X_i, X_j} \circ
c_{X_l, X_j} \circ c_{X_j, X_l} \circ c_{X_j, X_i} = \gamma_i\gamma_l \id_{X_j \ot X_i\ot X_l},
\]
where we omit identity morphisms and associativity constraints. Taking projection
of both sides on $X_k\ot X_j$ and computing the trace  we obtain  
the identity \eqref{ss=s}.

It is clear that the above maps between $G(\C)$ and
$\Aut_\otimes(\id_\C)$ are inverses of each other.
That they are group homomorphisms follows from the following
observation : if $X_j, X_{j_1}, X_{j_2}$
are invertible objects of $\C$ such that
$X_j \cong  X_{j_1} \otimes  X_{j_2}$ then
$s_{j_1 i} s_{j_2 i} = d_i s_{j i}$ for all $i\in \mathcal{I}$
(this is also a consequence of  the hexagon identity,
see \cite[Lemma 2.4(i)]{Mu}).
\end{proof}

\begin{remark}
For the representation category of a modular Hopf algebra
$H$ the result of Theorem~\ref{Schneider's iso} was proved
in \cite[Theorem 2.3(b)]{Sc}.
Namely, it was shown that the groups $G(H^*)$ and
$G(H)\cap Z(H)$ are isomorphic.
\end{remark}

\subsection{The lower and upper central series of a pseudounitary
modular category}

Let $\mathcal{K}$ be a full tensor subcategory of a braided tensor
category $\C$. In \cite{Mu} M.~M\"uger introduced the 
{\em centralizer} $\mathcal{K}'$ of $\mathcal{K}$, which is
a full tensor subcategory of $\C$. In the case when $\C$ is
pseudounitary and modular, $\mathcal{K}'$ is the fusion category 
generated by all simple objects
$X_i$ of $\C$  such that 
\begin{equation}
\label{sdd}
s_{ij} = d_i d_j
\end{equation}
for all simple  $X_j \in \mathcal{K}$.  If \eqref{sdd} holds
we will say that $X_i$ and $X_j$ {\em centralize} each other.


It was also shown in \cite{Mu} that  $\mathcal{K}'' =\mathcal{K}$
and  $\dim(\mathcal{K}) \dim(\mathcal{K}') =\dim(\C)$.
A fusion subcategory $\mathcal{K}\subseteq \C$
is symmetric if and only if
$\mathcal{K} \subseteq \mathcal{K}'$ and is modular
if and only if $\mathcal{K} \cap \mathcal{K}' = \Vec$.

 From now on we will assume that $\C$ is a pseudounitary modular category.
By \cite[Proposition 8.23]{ENO}, there is a canonical spherical structure
on $\C$ with respect to which the categorical dimensions of objects
coincide with their Frobenius-Perron dimensions: $d_i = \FPdim(X_i),
i\in \mathcal{I}$. In particular, each $d_i$ is a positive real number.

We will use the  notion of a centralizer
to relate the upper and lower
central series of a pseudounitary modular category.

\begin{lemma}
\label{inequality}
\begin{enumerate}
\item[(1)] For all $i,j\in \mathcal{I}$ we have $|s_{ij}| \leq d_i d_j$.
\item[(2)] The real part of $d_id_j - s_{ij}$ is non-negative and is equal
to $0$ if and only if $s_{ij} = d_id_j$.
\end{enumerate}
\end{lemma}
\begin{proof}
These results are essentially contained in \cite{Mu}; we give proofs
for the reader's convenience. The first inequality follows from  \eqref{sij}
and the Cauchy inequality:
\[
|s_{ij}| \ = | \theta_i^{-1}\theta_j^{-1} \sum_{k\in
\mathcal{I}}\, N_{i^*j}^k \theta_k d_k | \leq \sum_{k\in
\mathcal{I}}\, N_{i^*j}^k d_k = d_i d_j.
\]
Subtracting \eqref{sij} from $d_id_j = \sum_{k\in \mathcal{I}}\, N_{i^*j}^k d_k $ 
we obtain
\[
 d_id_j - s_{ij} =  \sum_{k\in \mathcal{I}}\, N_{i^*j}^k d_k (1- \theta_i^{-1} \theta_j^{-1} \theta_k).
\]
It is clear that the real part of every  summand in the last sum is non-negative. Hence,  the real part
of $ d_id_j - s_{ij}$ equals $0$ if and only if $\theta_k = \theta_i \theta_j$ for all $k\in\mathcal{I}$
such that $ N_{i^*j}^k \neq 0$. By \eqref{sij}, this is equivalent to $s_{ij} = d_id_j$.
\end{proof}

\begin{lemma}
\label{lemma1}
Let $X_i, X_k$ be simple objects in $\C$. Then
 $s_{kj} =d_kd_j$ for each $j\in\mathcal{I}$ such that $X_j$
is contained in $X_i\ot X_i^*$
if and only if $|s_{ik}|=d_id_k$.
\end{lemma}
\begin{proof}
Let us write the Verlinde formula \eqref{verlinde1}
in the following equivalent form
(cf.\ \cite[Lemma 2.4(iii)]{Mu}):
\begin{equation}
\label{Mugers eqn}
\frac{1}{d_k} s_{ki}s_{kl} = \sum_{j\in\mathcal{I}}\,
N_{il}^j s_{kj}\qquad \mbox{ for all } i,k,l\in \mathcal{I}.
\end{equation}
In particular, we have
\[
\frac{1}{d_k} |s_{ik}|^2  =
\frac{1}{d_k} s_{ik} s_{i^*k} = \sum_{j\in\mathcal{I}}\,
N_{ii^*}^j s_{kj}.
\]
If $s_{kj} =d_kd_j$ whenever $X_j$ is contained in $X_i\ot X_i^*$ then
\[
d_kd_i^2 = \sum_{j\in\mathcal{I}}\,
N_{ii^*}^j s_{kj} = \frac{1}{d_k} |s_{ik}|^2,
\]
whence $|s_{ik}| = d_id_k$.

Conversely, suppose $|s_{ik}| = d_id_k$. Then letting $l=i^*$ in
\eqref{Mugers eqn} we compute
\[
d_kd_i^2 = \sum_{j\in\mathcal{I}}\, N_{ii^*}^j s_{kj}.
\]
By Lemma~\ref{inequality}(2) the last equality is possible if and only if
$s_{kj}= d_kd_j$ for all $j$ such that $N_{ii^*}^j\neq 0$, so the 
proof is complete.
\end{proof}


\begin{proposition}
\label{key step}
Let $\mathcal{K}$ be a fusion subcategory of a pseudounitary modular
category $\C$. Let $\mathcal{L}$ be the fusion subcategory
generated
by simple objects $X_k\in \C$ such that $|s_{ik}| = d_i d_k$
whenever $X_i\in \mathcal{K}$.
Then $(\mathcal{K}_{ad})'= \mathcal{L} = (\mathcal{K}')^{co}$.
\end{proposition}
\begin{proof}
By definition, $(\mathcal{K}')^{co}$ is generated, as an abelian category,
by all simple $X_i\in \C$
such that every subobject of $X_i \otimes X_i^*$ centralizes $\mathcal{K}$.
By Lemma~\ref{lemma1}, the last condition is equivalent to 
$|s_{ik}| =d_id_k$ for all simple $X_k\in \mathcal{K}$.
Therefore, $(\mathcal{K}')^{co} = \mathcal{L}$.

Let $X_k$ be a simple object of $(\mathcal{K}_{ad})'$. Then
$X_k$ centralizes every simple subobject $X_j$ of $X_i \otimes X_i^*$
whenever $X_i$ is in $\mathcal{K}$. It follows from  Lemma~\ref{lemma1}
that $|s_{ki}| = d_k d_i$; i.e, $X_k \in \mathcal{L}$ and 
$(\mathcal{K}_{ad})' \subseteq  \mathcal{L}$.

Let $X_k$ be a simple object of $\mathcal{L}$, let $X_i$ be a simple
object of $\mathcal{K}$, and let $X_j$ be a simple subobject  of $X_i\otimes X_i^*$,
$k,j,i\in\mathcal{I}$. Then $s_{kj} = d_k d_j$ by  Lemma~\ref{lemma1}, and hence
$X_j$ belongs to $\mathcal{L}'$. This implies that $\mathcal{K}_{ad}$ is generated
by some of the simple objects of $\mathcal{L}'$, hence $\mathcal{K}_{ad} \subseteq  \mathcal{L}'$.

By M\"uger's double centralizer theorem \cite{Mu}, $\mathcal{L} = \mathcal{L}''$,
therefore  the above inclusions imply $\mathcal{L} = (\mathcal{K}_{ad})'$. 
\end{proof}

\begin{theorem}
\label{low = centralizer up}
Let $\C$ be a pseudounitary modular category and let
\begin{gather}
\label{up}
\C = \C^{(0)} \supseteq  \C^{(1)} \supseteq \cdots
\supseteq \C^{(n)} \supseteq \cdots \\
\label{lo}
\Vec = \C_{(0)} \subseteq  \C_{(1)} \subseteq \cdots
\subseteq \C_{(n)} \subseteq \cdots
\end{gather}
be the upper and lower central series of $\C$. Then
$(\C^{(n)})' = \C_{(n)}$ for all $n\geq 0$.
\end{theorem}
\begin{proof}
We have $(\C^{(0)})' =\C' =\Vec= \C_{(0)}$. If
$(\C^{(n)})' = \C_{(n)}$ for some $n\geq 0$ then,
using Proposition~\ref{key step} we obtain
$(\C^{(n+1)})' = ((\C^{(n)})_{ad})'=  ((\C^{(n)})')^{co} = ( \C_{(n)})^{co} =
\C_{(n+1)}$, and the statement follows by induction.
\end{proof}

Recall that $\C_{pt}$ denotes the
maximal pointed fusion subcategory of $\C$; i.e., the subcategory
generated by all invertible objects of $\C$.

\begin{corollary}\label{cor1}
Let $\mathcal{C}$ be a pseudounitary modular category. Then
$\C_{ad} = (\C_{pt})'$.
\end{corollary}
%

\begin{corollary}\label{cor2}
Let $\C$ be a nilpotent modular category of nilpotency class $c$.
Then for each $n\geq \frac{c}{2}$ the fusion subcategory
$\C^{(n)} \subseteq \C$ is symmetric.
\end{corollary}
\begin{proof}
Note that a nilpotent category is pseudounitary by Remark~\ref{long remark}(6).
Combining Theorems~\ref{upp=low} and \ref{low = centralizer up}
we obtain $\C^{(n)} = (\C_{(n)})'\subseteq (\C^{(c-n)})' \subseteq
(\C^{(n)})'$  for $n\geq \frac{c}{2}$.
\end{proof}

Recall that for any fusion category $\C$ its center $\mathcal{Z}(\C)$
is defined as the category whose objects are pairs $(X, c_{X, -})$ where
$X$ is an object of $\C$ and $c_{X,-}$ is a natural family of isomorphisms
$c_{X, V} : X\ot V \cong V \ot X$ for all objects $V$ in $\C$ satisfying 
certain compatibility conditions (see, e.g., \cite{K}). It is known that 
the center of a pseudounitary category is modular. 

Now assume that $\C$ is braided.
There are two full tensor embeddings of $\C$ into $\mathcal{Z}(\C)$ given
by $X\mapsto (X,\, c_{X,-})$ and by $X\mapsto (X,\, c^{-1}_{-,X})$. 
If  $\C_\pm$ denote the images of these embeddings then it is easy to see
that $\C_+' =\C_-$ in $\mathcal{Z}(\C)$.

\begin{theorem}\label{brnil}
Let $\C$ be a braided fusion category. Then $\C$ is nilpotent if and only
if its center $\mathcal{Z}(\C)$ is nilpotent. Moreover, if the nilpotency 
class of $\C$ is $c$ then the nilpotency class of $\mathcal{Z}(\C)$ is at 
most $2c$.
\end{theorem}
\begin{proof}
By Proposition~\ref{subcat+ surj image}, if $\mathcal{Z}(\C)$ is nilpotent, 
then so is $\C$.

Conversely, suppose $\C$ is nilpotent and identify it with 
$\C_+\subseteq \mathcal{Z}(\C)$. Let
\begin{equation}
\C= \C^{(0)} \supseteq  \C^{(1)} \supseteq \cdots
\supseteq \C^{(c)} =\Vec 
\end{equation}
be the upper central series of $\C$. Let $\E_k,\, k=0,\dots, c$, 
be the centralizer of $\C^{(k)}$ in $\mathcal{Z}(\C)$.  
By Proposition~\ref{key step}, we have 
\[
\E_k = (\C^{(k)})' =  ((\C^{(k-1)})_{ad})' =  ((\C^{(k-1)})')^{co} =(\E_{k-1})^{co}.
\]
Combining this with 
Lemma~\ref{coad in adco} we obtain $(\E_k)_{ad} = ((\E_{k-1})^{co})_{ad} \subseteq \E_{k-1}$
for all $k=0,\dots, c$. Note that $\E_c = (\C^{(c)})' = \mathcal{Z}(\C)$.
Assume that $\mathcal{Z}(\C)^{(k)} \subseteq \E_{c-k}$ for some $k\, (0\leq k < c)$.
Then $\mathcal{Z}(\C)^{(k+1)} =  (\mathcal{Z}(\C)^{(k)})_{ad} \subseteq (\E_{c-k})_{ad}
\subseteq \E_{c-k-1}$. Hence, $\mathcal{Z}(\C)^{(c)} \subseteq \E_0 = \C' = \C_-$ by induction.
Since $\C_-$ is nilpotent of class $c$, the result follows.
\end{proof}

\bibliographystyle{ams-alpha}

\end{document}